\documentclass[12pt]{amsart}
\usepackage{geometry}                
\usepackage{graphicx}
\usepackage{amssymb}
\usepackage{epstopdf}
 \usepackage{latexsym}
\usepackage{amsmath}
\usepackage{amsfonts} 
\usepackage{amsthm}

\newtheorem{Proposition}{Proposition}
\newtheorem{Theorem}[Proposition]{Theorem}
\newtheorem{Lemma}[Proposition]{Lemma}

\newtheorem{Remark}{Remark}

       \def\z{\noindent}

    \def\z{\noindent}  
 
    \def\sqr#1#2{{\vcenter{\vbox{\hrule height .#2pt
                             \hbox{\vrule width .#2pt height#1pt \kern#1pt
                                   \vrule width .#2pt}
                             \hrule height .#2pt}}}}

     \def\CC{\mathbb{C}}
    
    \def\NN{\mathbb{N}}

    \def\ZZ{\mathbb{Z}}

 \def\bfu{\mathbf{u}}
 \def\bfw{\mathbf{w}}
 \def\bff{\mathbf{f}}
  \def\bfg{\mathbf{g}}
 \def\bfn{\mathbf{n}}
 
  \def\bfh{\mathbf{h}}
    \def\bfP{\mathbf{P}}
 \def\bfp{\mathbf{p}}
   \def\bfm{\mathbf{m}}   
      \def\bfF{\mathbf{F}}
\def\bflam{\boldsymbol{\lambda}}
\def\bfphi{\boldsymbol{\phi}}
\def\bfpsi{\boldsymbol{\psi}}

\def\be{\begin{equation}}
\def\ee{\end{equation}}

\begin{document} 

\title[Nonlinear perturbations of Fuchsian systems]{Nonlinear perturbations of Fuchsian systems: corrections and linearization, normal forms}
\author{Rodica D. Costin}
        



\maketitle

\begin{abstract}Nonlinear perturbation of Fuchsian systems are studied in a region including two singularities. It is proved that such systems are generally not analytically equivalent to their linear part (they are not linearizable) and the obstructions are found constructively, as a countable set of numbers. Furthermore, assuming a polynomial character of the nonlinear part, it is shown that there exists a unique formal "correction" of the nonlinear part so that the "corrected" system is formally linearizable. 

Normal forms of these systems are found, providing also 
their classification.
\end{abstract}

\section{Introduction}\label{Introduction}

\subsection{Setting.}\label{Setting}
Nonlinear perturbations of Fuchsian systems are, in the present context, differential systems of the form
\be\label{genFuchSys}
\frac{d\bfu}{dx}=M(x)\bfu+{\bfg}(x,\bfu)
\ee
studied for $\bfu\in\CC^d$ small, and $x$ in a simply connected domain $D\subset\CC$ which includes singular points of the matrix $M(x)$.

It is assumed that the linear part of (\ref{genFuchSys})
\be\label{genLinFuchSys}
\frac{d\bfw}{dx}=M(x)\bfw\ \ \ \ \ \ \bfw\in\CC^d,\ x\in\CC
\ee
is a Fuchsian system. Then that all its singularities in the complex plane, including the point $x=\infty$, are of Fuchsian type (see the Appendix, \S\ref{RSP} for details). 

The function ${\bfg}(x,\bfu)$, which gathers the nonlinear terms, is assumed to be holomorphic, or, it can have (at most) simple poles at the singularities of $M(x)$: 
\be\label{defg}
{\bfg}(x,\bfu)=\frac{1}{Q(x)}{\bff}(x,\bfu)
\ee
where: $Q(x)=\prod(x-x_i)$, with $x_i\in D$ singular for $M$,  ${\bff}(x,.)$ has a zero of order two at $\bfu=0$  and  $\bff$ is holomorphic for $x\in D$ and small $\bfu$.

The case when the domain $D$ contains just one singular point of $M(x)$ was studied in \cite{Nonln} and the main results are summarized in \S\ref{rev1s}.

The present paper considers the case when $M(x)$ has three singularities in the extended complex plane. Placing one at $\infty$, the other two are conventionally placed at $x=1$ and $x=-1$ (but their location can be arbitrarily placed using a linear fractional change of the variable $x$). In this case the matrix $M$ has the form
\be\label{formM}
M(x)=\frac{1}{x-1}\, A\, +\, \frac{1}{x+1}\, B\ \ \ \ \ \  \ \,A,B\in\mathcal{M}_d(\CC)
\ee
and
\be\label{defQ}
Q(x)=x^2-1
\ee
(the condition that $x=\infty$ is also of Fuchsian type implies that the matrices $A$ and $B$ do not depend on $x$).

The systems studied have therefore the form
\be\label{perFuchs}
\frac{d\bfu}{dx}=\left(\frac{1}{x-1}\, A\, +\, \frac{1}{x+1}\, B\right)\bfu+\frac{1}{x^2-1}\bff(x,\bfu)
\ee
where $\bff$ is analytic for $x$ in the simply connected domain $D\ni\{\pm 1\}$ and for $\bfu\in\CC^d$, $|\bfu|<r$, $\bff(x,\cdot)$ has a zero of order two at $\bfu=0$.

In addition it is also assumed that the nonlinear part is of polynomial type in $x$, in the sense that
\be\label{pol_type}
{\mathbf{f}}(x,\bfu)=\sum_{|\mathbf{m}|\geq 2}\mathbf{f}_\mathbf{m}(x)\mathbf{u}^\mathbf{m}\ \ \ \ \ {\mbox{with}}\ \mathbf{f}_\mathbf{m}(x)\ {\mbox{polynomials}}
\ee

\subsection{Motivation.}\label{Motivation}

The question of linearization and, more generally, of equivalence,
 is a fundamental problem in the theory of differential equations. 
 
 Besides the clear theoretical interest, linearization and equivalence are used as instrumental methods in control theory, and in devising algorithms for numerical and symbolic calculations (see, to cite just a few authors, \cite{Fels}-\cite{Doubrov}). Often the methods used are developed from the method of equivalence introduced by Cartan \cite{Cartan} to decide whether two differential structures can be mapped one to another by a transformation taken in a given pseudogroup
\cite{Gardner}, \cite{Olver}. In the case of differential equations, the method was used for regular systems and extended for a neighborhood of one singular manifold \cite{Hermann}, \cite{Sternberg}.



The present paper uses direct and constructive methods in the study of linearizability and normal forms for equations in the class (\ref{perFuchs}) in domains which include two singular points.

The study of vector fields with an eigenvalue equal to $1$ at a singular point can be reduced to the study of a Fuchsian system near one singularity (\ref{gen1sing}) by eliminating time in a time-independent system.  Similarly, the study of vector fields in regions containing two, or more singular points is reducible to that of  equations of type (\ref{genFuchSys}).

More generally, the study of Hamiltonian systems with polynomial potentials near particular, periodic, or doubly-periodic solutions can be reduced to the study of nonlinear perturbations of Fuchsian systems \cite{MAA}, \cite{Thesis}.

The result of {Theorem}\,\ref{Obstructions} has a very interesting {{similarity with needed corrections found in other problems:}} \'Ecalle and Vallet showed that resonant systems are linearizable after appropriate correction \cite{Ecalle-Vallet}; also Gallavotti showed that there exists appropriate corrections of Hamiltonian systems so that the new system is integrable  \cite{Gallavotti}, convergence being proved later by Eliasson \cite{Eliasson}. This suggests the possible existence of an underlying general structure.

\subsection{Linearization in a neighborhood of one regular singular point and connection to integrability.}\label{Prior_results}

\subsubsection{Notations} For $\bfw=(w_1,\ldots ,w_d)\in\CC^d$ and $\bfm\in\NN^d$ denote
$$\bfw^\bfm=w_1^{m_1}\ldots w_d^{m_d}\ ,\ \ \ |\bfm|=m_1+\ldots +m_d\ ,\ \ \ \bfm\cdot\bfw=m_1w_1+\ldots +m_dw_d$$

\subsubsection{Region containing one regular singularity.}\label{rev1s}

Consider a nonlinear perturbation of a system with one regular singular point in the complex plane:
\be\label{gen1sing}
\frac{d\bfu}{dx}=\frac{1}{x}\, L(x)\bfu+\frac{1}{x}\, \tilde{\bff}(x,\bfu),\ \ \ \ \ \ \bfu\in\CC^d,\ x\in\CC
\ee
where ${L}(x)$ is a matrix analytic at $x=0$.

In nonresonant cases, after an analytic change of variables, it can be assumed that the matrix $L(x)$ is constant (see \S\ref{MxM0} for details):
\be\label{1sing}
\frac{d\bfu}{dx}=\frac{1}{x}\, L\bfu+\frac{1}{x}\, {\bff}(x,\bfu),\ \ \ \ \ \ \bfu\in\CC^d,\ x\in\CC
\ee
$\bff(x,\bfu)$ is assumed analytic for small $\bfu$: $|\bfu |<\rho$ and for $x$ in a domain $D_r$ which is either a disk centered at the origin: $|x|<r'$, or an annulus\footnote{Allowing the nonlinear part to be singular at $x=0$ accommodates systems corresponding to higher order equations.} $r''<|x|<r'$. 

Such systems are analytically linearizable if a Diophantine condition is satisfied \cite{Nonln}:

\begin{Theorem}\label{Th1sing}

Assume that the eigenvalues $\mu_1,...,\mu_d$ of the matrix $L$ satisfy the Diophantine
condition: there exist $C,\nu>0$ so that 
\be\label{DioCond}
\Big| \bfn\cdot\boldsymbol{\mu} +l-\mu_s\Big|>C\left(|\bfn|+|l|\right)^{-\nu}
\ee
for all $l\in\NN$ if $D_r$ is a disk (respectively all $l\in\ZZ$ if $D_r$ is an annulus), and for all $s\in\{1,...,d\}$, and 
$\mathbf{n}\in\NN^d$ with $ |\mathbf{n}|\geq 2$.

\z Then {{the system (\ref{1sing}) is analytically equivalent to its
linear part}}  $\mathbf{w}'=\frac{1}{x}L\mathbf{w}$ for $x\in D_{\tilde{r}}$ (a sub-disk, respectively, sub-annulus, of $D_r$: $0<\tilde{r}<r$) and $|\bfu |<\tilde{\rho}<\rho$.

\end{Theorem}

\begin{Remark}\label{remcirc}\ 

{{(i)}} The condition (\ref{DioCond}) is satisfied by generic matrices $L$ in the sense that the Lebesgue measure of the set of points $\boldsymbol{\mu}\in\CC^d$ which do not satisfy a condition (\ref{DioCond}) is zero, see \cite{Anosov-Arnold},Ch.6,\S\,3.

{{(ii)}} The domain $D_{\tilde{r}}$ can be made arbitrarily close to $D_r$ for $\tilde{\rho}$ small enough.

{{(iii)}} The analytic equivalence map with identity linear part is unique if no eigenvalue $\mu_j$ is integer.

\end{Remark}

\subsubsection{Region containing two singularities of Fuchsian type.}\label{inte}

Consider equation (\ref{perFuchs}). By Theorem \ref{Th1sing}, there (generically) exists an equivalence map of this equation to its linear part 
(\ref{genLinFuchSys}), for short, a {\em{linearization map}}, analytic at $\bfu=0$ and for $x$ in a neighborhood of $1$. This map is unique under the requirement of analyticity at $x=1$ (see Remark\,\ref{remcirc}). Similarly, there exists a unique linearization map which is analytic at $x=-1$. Generally these two maps do not coincide: the map which is analytic at $x=\pm 1$ is (usually) ramified at $x=\mp 1$. Indeed, the first term of the right-hand-side of (\ref{splithn}) below displays the ramification.

Therefore equations (\ref{genFuchSys}) are not necessarily linearizable in a region containing more than one singularity. 

Absence of linearizability has important consequences for the system, since it implies 
non-integrability (at least generically) \cite{RDC-MDK}:

\begin{Theorem}\label{RDC-MDK}
Consider the nonlinear equation in one dimension:
\be\label{rdc_mdk_art}
{\frac{du}{dx}=\left(\frac{a}{x-1}+\frac{b}{x+1}\right)u+\frac{1}{x^2-1}\, {f(x,u)}}
\ee
for $x$ in a domain $D$ containing $x=\pm1$ and small $u$.

{{If}} equation (\ref{rdc_mdk_art})
 {is {\em{not}} analytically linearizable}
 then for generic $a,b$ (precise conditions are given in \cite{RDC-MDK})
solutions have dense branching and therefore 
{no single-valued integrals exist.}
 Among integrable cases, first integrals are not meromorphic {{(generically)}}.

\end{Theorem}

It is therefore important to find what are the linearization criteria for equations (\ref{genFuchSys}), and furthermore to describe the equivalence classes - to find normal forms for these equations.

\section{Main Results}\label{Main_res}

\subsection{Assumptions.}\label{eq}

Consider the system (\ref{perFuchs}). The following conditions (a)-(c), satisfied by generic matrices, are assumed.

(a) The eigenvalues of $A$, respectively $B$, satisfy the Diophantine condition 
(\ref{DioCond}), and 
(b) none of these eigenvalues is an integer.

(c) The eigenvalues $\lambda_1,\ldots ,\lambda_d$ of the matrix $A+B$ satisfy the following nonresonance condition:
\be\label{non_res}
k+\mathbf{n}\cdot {\bflam}-\lambda_j \not=0\ \ \ {\mbox{for\ all}}\ \ \mathbf{n}\in\NN^d\ ,\ k\in\NN\ ,\ \ j=1,\ldots ,d
\ee

\subsection{Obstructions to linearization and existence of corrections.}\label{SecObs}

\begin{Theorem}\label{Obstructions}

Consider the system (\ref{perFuchs}), (\ref{pol_type}) under the assumptions of \S\ref{eq}.

Then {\em{there exists a unique}} "correction" $\bfphi(\bfu)$ of $\bff(x,\bfu)$ as a formal series 
\begin{equation}\label{serphi1}
\bfphi(\bfu)=\sum_{|\bfm|\geq 2}\,  \bfphi_\bfm \bfu^\bfm
\end{equation}
so that the "corrected" equation 
\be\label{corperFuchs}
\frac{d\bfu}{dx}=\left(\frac{1}{x-1}\, A\, +\, \frac{1}{x+1}\, B\right)\bfu+\frac{1}{x^2-1}\left[\bff(x,\bfu)-\bfphi(\bfu)\right]
\ee
is linearizable by a formal series
\be\label{serh}
{\bfu=\mathbf{H}(x,\bfw)=\bfw+\sum_{|\mathbf{m}|\geq 2}\mathbf{h}_\bfm(x)\mathbf{w}^\mathbf{m}}
\ee
where $\mathbf{h}_\bfm(x)$ are functions analytic on $D$. 

Moreover,  {$\mathbf{h}_\bfm(x)$ are polynomials}.

\end{Theorem}

{\bf{Remarks}} 

{\bf{1.}} The coefficients $\bfphi_\mathbf{m}\in\CC^d$ of the correction $\bfphi$ (see (\ref{serphi1})) represent the obstructions to linearization, in the sense that the system (\ref{perFuchs}) is linearizable by a formal power series if and only if all these vectors are zero.

Obviously, if no formal power series (\ref{serh}) exists, then analytic equivalence is precluded too.

{\bf{2.}}  $\bfphi_\mathbf{m}$ are obtained constructively, see (\ref{eqhn}), (\ref{formRn}), Lemmas\,\ref{FundamLemma} and\,\ref{sollem}.

{\bf{3.}} If the nonresonance condition (\ref{non_res}) does not hold then existence, or uniqueness, of corrections may fail, as illustrated in the Appendix \S\ref{reson}. 

The condition that $\bff(\cdot,\bfu)$ be of polynomial type is expected to generalizable to analytic functions. In a forthcoming paper \cite{Analytic} it is shown that if the matrices $A$ and $B$ are simultaneously diagonalizable, then {Theorem}\,\ref{Obstructions} holds more generally, for $\bff(x,\bfu)$ analytic in $x$.
Furthermore, it is also shown in \cite{Analytic} that if in addition, the eigenvalues of the matrices $A$ and $B$ have positive real part, then all the series converge and the corrected system is analytically linearizable.


\subsection{Formal normal forms.}

Since equations (\ref{perFuchs}) are not necessarily linearizable, then two such equations are not necessarily equivalent. The following theorem finds normal forms for equations (\ref{perFuchs}).

\begin{Theorem}\label{NormForm}

Consider the system (\ref{perFuchs}), (\ref{pol_type}) under the assumptions of \S\ref{eq}.

Then {\em{there exists a unique}} formal series $\bfpsi(\bfw)=\sum_{|\bfm|\geq 2}\bfpsi_\bfm\bfw^\bfm$ so that  (\ref{perFuchs}) is equivalent to
\be\label{NF}
\bfw'=\left(\frac{1}{x-1}A+\frac{1}{x+1}B\right)\bfw+\frac{1}{x^2-1}\bfpsi(\bfw)
\ee

\z through a formal series
\be\label{hNF}
\bfu=\mathbf{H}(x,\bfw)=\bfw+\sum_{|\bfm|\geq 2}\mathbf{h}_\bfm(x)\bfw^\bfm
\ee
 with 
$\mathbf{h}_\bfm(x)$ analytic functions on $D$. 

Moreover,  $\mathbf{h}_\bfm(x)$ are polynomials.
\end{Theorem}

\section{Proofs}\label{Proofs}

\subsection{Proof of Theorem \ref{Obstructions}.}\label{PfOb}

\subsubsection{The recursive system.}\label{recsys}

A map of the form (\ref{serh}), $\bfu=\mathbf{H}(x,\bfw)=\bfw+\mathbf{h}(x,\bfw)$, is a linearization map of (\ref{perFuchs}) if and only if $\bfh$ satisfies the nonlinear partial differential equation
\be\label{PDE}
\partial_x\bfh+d_\bfw\bfh\, M\bfw=M\bfh+\frac{1}{x^2-1}\, \left[\bff(x,\bfw+\bfh)-\bfphi(\bfw+\bfh)\right]
\ee
where $M$ is given by (\ref{formM}).

Denote by $\bfh_n$ the homogeneous part of degree $n$ of a function $\bfh(x,\bfw)$:
\be
\bfh_n(x,\bfw)=\sum_{|\mathbf{m}|=n}\bfh_\mathbf{m}(x)\bfw^\mathbf{m}
\ee
(and a similar notation is used for $\bff(x,\bfw)$ etc.).

Writing solutions of (\ref{PDE}) as power series in $\bfw$ in  we obtain $\bfh_n$ recursively in $n$. The functions $\{\bfh_\bfm\}_{|\bfm|=n}$ satisfy the linear differential system 
\be\label{eqhn}
\partial_x\bfh_n+d_\bfw\bfh_n\, M\bfw-M\bfh_n=\frac{1}{x^2-1}\, \mathbf{R}_n(x,\bfw)
\ee
where
\be\label{formRn}
\mathbf{R}_n=\bff_n-\bfphi_n+\tilde{\mathbf{R}}_n
\ee
with $\tilde{\mathbf{R}}_n$ a
polynomial in $\bfphi_\mathbf{m}$, $\bfh_\bfm$, $\bff_\mathbf{m}$ with $|\mathbf{m}|<n$.

\subsubsection{Numerical obstructions to linearizability: an analytic illustration}\label{numob}

The existence of the correction $\boldsymbol{\phi}(\bfu)$ is established below in \S\ref{extphi} using algebraic methods which rely on the polynomiality assumption (\ref{pol_type}). It is interesting, however, to look at the analytic context, as a hint on why existence, and uniqueness, of this correction is to be expected.

Denote by $Y(x)$ a fundamental matrix of the linear system (\ref{genLinFuchSys}): $Y'=MY$. If $Q^{-1}Y^{-1}$ is Lebesgue integrable at $x=\pm 1$, then 
a particular solution of (\ref{eqhn}) is given by
\be\label{hnat-1}
\bfh_n(x,\bfw)=Y(x)\int_{-1}^xQ(t)^{-1}Y(t)^{-1}\mathbf{R}_n(t,Y(t)Y(x)^{-1}\bfw)\, dt
 \ee

Let {$G$} be the {monodromy matrix} of (\ref{genLinFuchSys}) at $x=-1$: after analytic continuation along a closed loop around $x=-1$ the matrix $Y(x)$ becomes {$AC_{(-1)}Y(x)=Y(x)G$}.

The analytic continuation of (\ref{hnat-1}) on a closed loop around $x=-1$ yields
\be
AC_{(-1)}{\bfh_n(x,\bfw)=Y(x)G\int_{-1}^xQ(t)^{-1}G^{-1}Y(t)^{-1}\mathbf{R}_n(t,Y(t)GG^{-1}Y(x)^{-1}\bfw)dt}
\ee
\be\nonumber
=\bfh_n(x,\bfw)
\ee

\z which means that (\ref{hnat-1}) are the homogeneous terms of the unique linearization map of  (\ref{perFuchs}) which is analytic at $x=-1$. 

Rewriting (\ref{hnat-1}) as
\be\label{splithn}
\bfh_n(x,\bfw)=Y(x)\int_{-1}^1Q^{-1}Y^{-1}\mathbf{R}_ndt+Y(x)\int_{1}^xQ^{-1}Y^{-1}\mathbf{R}_ndt
 \ee
a similar argument shows that the last term in (\ref{splithn}) is the homogeneous term of the unique linearization map which is analytic at $x=1$. Then $\bfu=\bfw+\mathbf{h}(x,\bfw)$ is analytic at both $x=1$ and $x=-1$ if and only if
\be\label{obcond}
\int_{-1}^1Q(t)^{-1}Y(t)^{-1}\mathbf{R}_n(t,Y(t)Y(x)^{-1}\bfw)\,  dt\, =\, 0\ \ \ {\mbox{for\ all\ }}\bfw\in\CC^d,\ n\geq 2
\ee

The recursive formulas (\ref{obcond}) contain the obstructions to linearization: there are $d$ numerical conditions (one vector in $\CC^d$) for every $\mathbf{m}\in\NN^d, |\mathbf{m}|\geq 2$.

\subsubsection{Existence of the correction $\boldsymbol{\phi}$.}\label{extphi}

The existence of $\boldsymbol{\phi}$ (as a formal series) is established directly on the differential equation  (\ref{eqhn}) with (\ref{formRn}), rather than solving (\ref{obcond}). This result, which completes the proof of Theorem\,\ref{Obstructions}, follows directly from the following lemma:

\begin{Lemma}\label{FundamLemma}

Denote by $\mathcal{P}_n$ the space of vector-valued polynomials in $\bfw$, homogeneous degree 
$n$ ($\mathcal{P}_n\subset\CC^d[w_1,\ldots ,w_d]$).

Assume (\ref{non_res}) holds.

Then for any $\bfF=\bfF(x,\bfw)\in \mathcal{P}_n[x]$, there exists a unique $\bfphi\in\mathcal{P}_n$ so that the differential equation 
\be\label{fdlem}
 \partial_x\bfP(x,\bfw)+d_\bfw\bfP(x,\bfw)\, M(x)\bfw-M(x)\bfP(x,\bfw)=\frac{1}{x^2-1}\, \left[\bfF(x,\bfw)-\bfphi(\bfw)\right]
\ee
has a solution $\bfP\in\mathcal{P}_n[x]$.

\end{Lemma}

The proof of Lemma \ref{FundamLemma} is given in \S\ref{pfLemma}.

Note that if $\bfF$ has degree ${k+1}$ in $x$, then the solution $\bfP$ of (\ref{fdlem}) has degree $k$. 

\begin{Remark}\label{Rodrigues}

The only solutions of (\ref{fdlem}) which are analytic at both $x=1$ and $x=-1$ are polynomials, and they are given by the following generalized Rodrigues formulas. Consider the polynomials $\bfF_{k+1}$ given by
$$\bfF_{k+1} \left(x,\bfw;\mathbf{q}\right) = Q(x)Y(x) \frac{d^{k+1}}{dt^{k+1}}\Big|_{t=x}\  \left[ Q(t)^kY(t)^{-1}\mathbf{q}(Y(t)Y(x)^{-1}\bfw) \right] \ ,\  \ \ \ \mathbf{q}\in\mathcal{P}_n\ ,\ k\in\NN$$
(these, together with $\mathcal{P}_n$,  generate the vector space $\mathcal{P}_n[x]$). Then 
\be\label{PJ}
\bfP_k \left(x,\bfw;\mathbf{q}\right) = Y(x) \frac{d^{k}}{dt^{k}}\Big|_{t=x}\  \left[ Q(t)^kY(t)^{-1}\mathbf{q}(Y(t)Y(x)^{-1}\mathbf{w}) \right]
\ee
is the polynomial solution of (\ref{fdlem}) for $\bfF=\bfF_{k+1}$ and $\bfphi=0$. 

The functions (\ref{PJ}) are indeed polynomials in $x$, of degree $k$. For example
$$\bfP_0 \left(x,\bfw;\mathbf{q}\right) = \mathbf{q}(\bfw),\ \ \bfP_1 \left(x,\bfw;\mathbf{q}\right) = (Q'-QM)\,\mathbf{q}(\bfw)+Q\, d\mathbf{q}(\bfw)\, M\bfw$$
$$\bfP_2 \left(x,\bfw;\mathbf{q}\right) =\left[ (Q^2)''+Q^2(M^2-M')-2(Q^2)'M\right] \mathbf{q}(\bfw)+2\left[(Q^2)'+Q^2M\right] \, d\mathbf{q}(\bfw)\, M\bfw$$
$$+Q^2\, d\mathbf{q}(\bfw)\, (M'+M^2)\bfw+Q^2\, d^2\mathbf{q}(\bfw)\, (M\bfw,M\bfw)$$

\end{Remark}


\subsubsection{Proof of Lemma \ref{FundamLemma}.}\label{pfLemma}

The proof uses the following solvability Lemma:

\begin{Lemma}\label{sollem}

Let $\Lambda$ be a $d\times d$ matrix with eigenvalues satisfying the nonresonance condition 
(\ref{non_res}). Let $J_\Lambda$ be the following linear operator on $\mathcal{P}_n$:
\be\label{defJ}
(J_\Lambda\bfp)(\bfw)=(d\bfp )\, \Lambda \bfw-\Lambda\bfp
\ee

Let $l\in\NN$.

(i) If $\Lambda$ is diagonal, then the operator $l+J_\Lambda$ is diagonal, with eigenvectors $\bfw^\bfm{\bf{\rm{e}}}_j$, and corresponding eigenvalues $l+\bfm\cdot\boldsymbol{\lambda}-\lambda_j$, for all $j=1,\ldots,d$ and $\bfm$, $|\bfm|=n$.

(ii)  If $\Lambda$ is in Jordan normal form, then $l+J_\Lambda$ is in Jordan normal form and the eigenvalues are as at (i).

(ii) As a consequence, the operator $l+J_\Lambda$ is one-to-one and onto $\mathcal{P}_n$.

\end{Lemma}

Lemma \ref{sollem} follows by calculation (see also \cite{Arnold} Ch.5, \S{22}).

To prove Lemma \ref{FundamLemma} denote by $L$ and $N$ the following linear operators on $\mathcal{P}_n$: $L=J_{(A+B)}$ and $N=J_{(A-B)}$ (see (\ref{defJ})).

Let $\bfF=\bfF(x,\bfw)\in \mathcal{P}_n[x]$ be a polynomial in $x$ degree $k+1$:
$\bfF=\sum_{j=0}^{k+1}x^j\bff_j$ with $\bff_j\in\mathcal{P}_n$.
Identifying the coefficients of $x$, a solution $\bfP$  of (\ref{fdlem}) has degree $k$ and, if $\bfP=\sum_{j=0}^{k}x^j\bfp_j$  then $\bfp_j$ must satisfy the recursive system
\be\label{sysP}
\begin{array}{ll}
(k+L)\bfp_k=\bff_{k+1}\\
(k-1+L)\bfp_{k-1}=N\bfp_k+\bff_{k}\\
(l-1+L)\bfp_{l-1}=N\bfp_l+(l+1)\bfp_{l+1}+\bff_{l}\ \ \ {\mbox{for\ }} 1\leq l\leq k-1\\
0=N\bfp_0+\bfp_1+\bff_0-\bfphi
\end{array}
\ee

By Lemma \ref{sollem} the system (\ref{sysP}) can be solved recursively for $\bfp_k$, then $\bfp_{k-1}$, and all the way to $\bfp_0$, for any $\bff_{k+1}$, $\bff_k$, $\ldots$, $\bff_1$. The last relation in (\ref{sysP}) determines uniquely $\bfphi$. \qed

\subsection{Proof of {Theorem} \ref{NormForm}.}

The map (\ref{hNF}) is an equivalence map of (\ref{perFuchs}) and (\ref{NF}) iff $\bfh$ satisfies the nonlinear PDE
\be
\partial_x\bfh+d_\bfw\bfh\, M\bfw-M\bfh=\frac{1}{x^2-1}\, \left[\bff(x,\bfw+\bfh)-\bfpsi(\bfw)-\left(d_\bfw\bfh\right)\bfpsi\right]
\ee
where $M$ is given by (\ref{formM}).
As in the proof of Theorem \ref{Obstructions}, expanding in power series in $\bfw$
and denoting by $\bfh_n$ the homogeneous part of degree $n$ in (\ref{serh}) 
we obtain $\bfh_n$ recursively as solutions of the equations
\be\label{eqhnp}
\partial_x\bfh_n+d_\bfw\bfh_n\, M\bfw-M\bfh_n=\frac{1}{x^2-1}\, \mathbf{R}_n(x,\bfw)\ \ \ \ \ (n\geq 2)
\ee
where
\be\label{formRnp}
\mathbf{R}_n=\bff_n-\bfpsi_n+\tilde{\mathbf{R}}_n
\ee
with $\tilde{\mathbf{R}}_n$ a
polynomial in $\bfpsi_\mathbf{m}$, $\bfh_\bfm$, $\bff_\mathbf{m}$ with $|\mathbf{m}|<n$.

The terms $\bfpsi_n$ are  recursively determined from (\ref{eqhnp}), (\ref{formRnp}) using  Lemma \ref{FundamLemma} for $\bfF=\bff_n+\tilde{\mathbf{R}}_n$.\qed


\section{Acknowledgments}

The author is grateful to the two anonymous referees for the careful reading of the manuscript and for their comments and suggestions.

\section{Appendices}

\subsection{Singular points of Fuchsian type of a linear differential equation.}\label{RSP}

Consider a linear differential equation 
\be\label{olineq}
\frac{d\bfw}{dx}=M(x)\bfw\ \ \ \ \ \ \bfw\in\CC^d,\ x\in\CC
\ee

A point $x=x_0$ is called a regular point of (\ref{olineq}) if the matrix $M(x)$ is analytic at $x=x_0$. 

A point $x=x_0$ is called {{singular of Fuchsian type}} if $M(x)$ has an isolated singularity at $x=x_0$ and there exists a constant $k$ so that all solutions of the system, in every sector of the complex $x$-plane with vertex at $x_0$ grow no faster than $x^k$ as $x\to x_0$ (see \cite{Anosov-Arnold}, Ch.7, \S 2.1). 

The equation (\ref{olineq}) is called Fuchsian if 
$$M(x)=\sum_{j=1}^r\frac{1}{x-a_j}A_j$$
in which case solutions of (\ref{olineq}) are analytic at all points except $a_1,\ldots ,a_r$ and $\infty$ and these points are of Fuchsian type.
It turns out that in this case all solutions have formal asymptotic expansions as series in powers and logs, as $x\to a_j$, or $x\to\infty$  \cite{Wasov}. Moreover, these series converge.

For {\em{generic}} nonlinear perturbations of systems (\ref{olineq}) near a singular point of Fuchsian type $x_0$, solutions also have convergent asymptotic expansions in powers of $x$ as $x\to x_0$ by Theorem\,\ref{Th1sing}.
For simplicity, the point $x=x_0$ will be also called a singular point of Fuchsian type (of a nonlinear system).

\subsection{The matrix $L(x)$ can be assumed constant in (\ref{gen1sing}).}\label{MxM0}

Indeed, if the eigenvalues of $L(0)$ are nonresonant (in the sense that no two eigenvalues differ by an integer) then the linear part of (\ref{gen1sing}):
\be\label{Mx}
 \frac{d{\bfw}}{dx}=\frac{1}{x}\, L(x){\bfw}
\ee
is analytically equivalent to 
\be\label{M0}
 \frac{d\tilde{\bfw}}{dx}=\frac{1}{x}\, L(0)\tilde{\bfw}
\ee
for small $x$ and $\bfw$ \cite{Anosov-Arnold}. Denoting ${\bfw}=H(x)\tilde{\bfw}$, the analytic equivalence map between (\ref{Mx}) and (\ref{M0}), the same map transforms (\ref{gen1sing}) into a system with the same regularity, but with $L$ constant:  the map ${\bfu}=H(x)\tilde{\bfu} $ is an analytic equivalence of (\ref{gen1sing}) and
\be\label{M01sing}
\frac{d\tilde{\bfu}}{dx}=\frac{1}{x}\, L(0)\tilde{\bfu}+\frac{1}{x}\, H^{-1}\tilde{\bff}(x,H\tilde{\bfu})
\ee

\subsection{Examples in which condition (\ref{non_res}) does not hold.}\label{reson}

To illustrate that existence of correction, or their uniqueness, may fail if the condition (\ref{non_res}) is not satisfied consider the simplest nonlinear equations: Riccati equations (in dimension one)
\be\label{Riccati}
u'(x)\, =\, \left( \frac{a}{x-1}+\frac{b}{x+1}\right)u(x)\, +\, \frac{f(x)}{x^2-1}\, u^2(x)
\ee

For $n=2$ the equation (\ref{eqhn}) is
\be\label{oeq}
{h_2'(x)}\, +\, \left( \frac{a}{x-1}+\frac{b}{x+1}\right)\, h_2(x)\, =\, \frac{f(x)-\phi_2}{x^2-1}
\ee

Take for simplicity $f(x)=mx^2+nx$ and $a=b$ and note that for $a=b=-1/2$ the condition (\ref{non_res}) is not satisfied. 
A simple calculation yields that equation (\ref{oeq}) has an analytic solution for $\phi_2=m/(2a+1)$ and this solution is $h_2(x)=m/(2a+1)x+n/(2a)$. Furthermore $\phi_n=0$ for $n>2$ and it can be easily checked that $h_n=h_2^{n-1}$.

These formulas certainly cannot be used in the resonant case $a=b=-1/2$. In this case the general solution of (\ref{oeq}) is
\be\label{11}
h_2(x)=m\left[ (1-x^2)^{1/2}\arcsin x-x\right]-n+\phi_2x+C(1-x^2)^{1/2}
\ee

If $m\not=0$ then no solution (\ref{11}) is single-valued at both $x=\pm 1$, no matter what $\phi_2$ is: existence of corrections fails.

By contrast, for $m=0$ uniqueness fails, since the solution (\ref{11}) with $C=0$ is analytic for any $\phi_2$.


\end{document}